\documentclass[final]{siamltex}

\usepackage{amsmath}
\usepackage{mathrsfs}
\usepackage{amsfonts}
\usepackage{amsmath,amsfonts,amsbsy,amstext,amsopn,amssymb}

\DeclareMathOperator{\e}{e}
\DeclareMathOperator{\B}{B}

\def\Expect{\mathop{\mathbf E}}
\def\Prob{\mathop{\mathrm{Prob}}}

\DeclareMathOperator{\Vol}{Vol}

\newcommand{\I}{I}

\newcommand{\R}{\mathbb{R}}

\newcommand{\PP}{\mathbb{P}}

\newcommand{\Ccal}{\mathscr C}
\newcommand{\Ocal}{\mathscr O}
\newcommand{\Id}{\mathsf{Id}}
\def\BAMS{Bulletin of the Amer. Math. Soc.}

\def\e{\varepsilon}

\edef\qedrestoreat{\noexpand\catcode\lq\noexpand\@=\the\catcode\lq\@}
\catcode\lq\@=11

\ifx\protect\undefined\let\protect\relax\fi

\def\qed{\protect\@qed{$\qedsymbol$}}
\def\pushright{\protect\@pushright}

\def\QED{\protect\@qed{{\rm Q.E.D.}}}

\def\QEI{\protect\@qed{{\rm Q.E.I.}}}

\def\Proof{\protect\@Proof}\def\endProof{\protect\@endProof}%

\def\Proofof#1{\protect\@Proofof{#1}}\def\endProofof{\protect\@endProofof}%

\def\qedsymbol{\raisebox{-.2ex}{$\Box$}}

\def\TheWordProof{\sc Proof.}
\def\TheWordProofof#1{\sc Proof of #1.}

\def\ProofFont{}

\newif\ifAutoQED\AutoQEDfalse

\newif\ifNumberResults
\ifx\theorem@style\undefined
   \NumberResultstrue
   
\fi

\def\parag@pushright#1{{
    \parfillskip=0pt            
    \widowpenalty=10000         
    \displaywidowpenalty=10000  
    \finalhyphendemerits=0      
    \hbox@pushright             
    #1
    \par}}

\def\hbox@pushright{
    \unskip                     
    \nobreak                    
    \hfil                       
    \penalty50                  
    \hskip.2em                  
    \null                       
    \hfill                      
}%

\newif\if@qed\@qedfalse
\def\save@set@qed{\let\saved@ifqed\if@qed\global\@qedtrue}%
\def\restore@qed{\global\let\if@qed\saved@ifqed}

\def\@Proof{%
   \par\removelastskip\bigskip\penalty100
   \save@set@qed
   \noindent\ProofFont{\TheWordProof\enskip}%
}%
\def\@Proofof#1{%
   \par\removelastskip\bigskip\penalty100
   \save@set@qed
   \noindent\ProofFont{\TheWordProofof{#1}\enskip}%
}%
\def\@endProof{%
   \qed\restore@qed
   \penalty-100 \medskip
}
\def\@endProofof{%
   \qed\restore@qed
   \penalty-100 \medskip
}
\def\@qed#1{%
\if@qed                                 
     \global\@qedfalse
        \ifmmode\ifinner\pushright{#1}
        \else\eqno{\qedsymbol}\fi
        \else\pushright{#1}\fi%
\else\ifhmode\ifinner\else\par\fi\fi
\fi}
\def\@pushright#1{%
  {\ifvmode                       
       \null\hfill{#1}\par        
  \else\ifmmode\maths@pushright{\hbox{#1}}
       \else\ifinner\hbox@pushright{#1}
            \else\parag@pushright{#1}
  \fi  \fi  \fi
}}%
\def\maths@pushright#1{{%
  \ifinner
     \hbox@pushright{#1}%
  \else
     \eqno#1
     \def\]{$$\ignorespaces}
  \fi
}}%
\begin{document}

\title{Adversarial Smoothed Analysis}

\author{Felipe Cucker\footnotemark[2],
Raphael Hauser\footnotemark[1] and Martin Lotz\footnotemark[1]}

\renewcommand{\thefootnote}{\fnsymbol{footnote}}
\footnotetext[1]{Oxford University Computing Laboratory, Wolfson
Building, Parks Road, Oxford, OX1 3QD, United Kingdom,
(hauser@comlab.ox.ac.uk).}
\footnotetext[2]{City University of Hong Kong,
Department of Mathematics, Kowloon Tong, HONG KONG,
(macucker@cityu.edu.hk). Partially supported by GRF grant
CityU 100808.}
\renewcommand{\thefootnote}{\arabic{footnote}}

\date{\today}
\maketitle  

\begin{abstract}
The purpose of this note is to extend the results on uniform
smoothed analysis of condition numbers from \cite{BuCuLo:07} to
the case where the perturbation follows a radially symmetric
probability distribution. In particular, we will show that the
bounds derived in~\cite{BuCuLo:07} still hold in the case of
distributions whose density has a singularity at the center of the
perturbation, which we call {\em adversarial}.
\end{abstract}

\begin{AMS}
Primary 65Y20; secondary 65G99.
\end{AMS}
\begin{keywords} Condition numbers, random matrices,
average case analysis, smoothed analysis.
\end{keywords}

\section{Introduction}\label{se:introduction}

Condition numbers play a central role in numerical analysis. They
occur in error analysis for finite-precision algorithms (this
being historically the reason for their introduction in the late
1940's by von Neumann and Goldstine~\cite{vNGo47} and
Turing~\cite{Turing48}) as well as a parameter in expressions
bounding the number of iterations in a variety of algorithms (a
paradigmatic example being the conjugate gradient
method~\cite[Theorem~38.5]{TrefethenBau}). In practice, however, a
difficulty appears: it would seem that to know the condition
number of a given data one needs to solve the problem at hand on
this data. An inconvenient circularity. A way out of it, proposed
by Steve Smale (see~\cite{Smale97} for a review), is to assume a
probability measure on the space of data and to study the
condition number $\Ccal(a)$ at data $a$ as a random variable. In
other words, to study the condition number of random data.

In doing so Demmel~\cite{Demmel87} noticed that most
condition numbers could be written as (or at least reasonably
sharply bounded by) the relativized inverse of the distance from
the data $a\in\R^{n+1}$ to a set of ill-posed instances
$\Sigma\subset\R^{n+1}$.  That is, one
could write
\begin{equation}\label{eq:conic}
   \Ccal(a)=\frac{\|a\|}{\mathrm{dist}(a,\Sigma)}.
\end{equation}
The simplest example of this phenomenon is given by the condition
number for matrix inversion and linear equation solving. For a
non-singular square matrix $A$ it takes the form $\kappa
(A):=\|A\|\|A^{-1}\|$, where $\| \; \|$ denotes the operator norm.
The Condition Number Theorem by Eckart and Young states that
$\|A^{-1}\|=d(A,\Sigma)^{-1}$, where $\Sigma$ is the set of
singular matrices.

In most applications, $\Sigma$ is a pointed cone. Therefore, one
could normalize so that $a$ belongs to the $n$-dimensional unit
sphere $S^n$. Note that the usual assumption that $a$ has a
Gaussian distribution in $\R^{n+1}$ yields a uniform distribution
in $S^n$ after this normalization. It is for condition numbers as
in~\eqref{eq:conic}
---which we shall call {\em conic}--- with inputs drawn from the uniform
distribution on $S^n$ that Demmel proved in~\cite{Demmel88}
(shortly after~\cite{Demmel87}) a general result bounding their
tail as a function of $n$ and the degree of an algebraic hypersurface
containing $\Sigma$.

Very recently, a new paradigm for probabilistic analysis was
proposed by Spielman and Teng~\cite{ST:02,ST:06}. Called {\em
smoothed analysis}, it consists of replacing the idea of ``random
data'' by that of ``random perturbation of a given data'' and
study the worst-case (w.r.t.\ data $a$) of the latter. In its
original formulation, and in the case of a condition number
$\Ccal(a)$, this amounts to study the tail
\begin{equation*}
     \sup_{a\in\R^{n+1}} \Prob_{z\in N(a,\sigma^2)}
           \left\{\Ccal(z)\geq t\right\}
\end{equation*}
or the expected value
\begin{equation*}
\sup_{a\in \R^{n+1}}\,\Expect_{z\in N(a,\sigma^2)}[\ln\Ccal(z)]
\end{equation*}
where $N(a,\sigma^2)$ is a Gaussian distribution centered at $a$
with covariance matrix $\sigma^2\Id$ and $\sigma^2$ small (with
respect to $\|a\|$). In~\cite{BuCuLo:07}, to obtain general
results as in~\cite{Demmel88}, data was again restricted to $S^n$
and the expressions above replaced by
\begin{equation*}
     \sup_{a\in S^n} \Prob_{z\in B(a,\sigma)}
           \left\{\Ccal(z)\geq t\right\}
\end{equation*}
and
\begin{equation*}
\sup_{a\in S^{n}}\,\Expect_{z\in B(a,\sigma)}[\ln\Ccal(z)]
\end{equation*}
where $B(a,\sigma)$ is the open ball (that is, the spherical cap)
in $S^n$ centered at $a$ and of radius $\sigma$, and $z$ is drawn
from a uniform distribution on this ball.

One of the claimed advantages of smoothed analysis is a smaller
dependence on the underlying distribution. It follows from this
claim that the replacement of Gaussian perturbations by uniform
ones should not significantly affect the smoothed analysis of
$\Ccal(a)$. The goal of this note is to further pursue this claim
by extending the main result in~\cite{BuCuLo:07}, combining it
with ideas from~\cite{HM:06}, to a class of distributions we call
{\em adversarial}. The support of such a distribution is, as in
the uniform case, the ball $B(a,\sigma)$ and they are radially
symmetric as well. But their density increases when approaching
$a$ and has a pole at $a$.

\section{Preliminaries}\label{se:preliminaries}
We assume our data space is $\R^{n+1}$, endowed with a scalar
product $\langle \ ,\ \rangle$. In all that follows we consider
problems whose set of ill-posed inputs $\Sigma$ is a point-symmetric
cone in $\R^{n+1}$. That is, if $x\in\Sigma$ then $\lambda
x\in\Sigma$ for all $\lambda\in\R$. By a {\em conic condition
number} we understand a function $\Ccal\colon \R^{n+1}\rightarrow
[1,\infty]$ such that for all $a\in \R^{n+1}$ we have
\begin{equation*}
  \Ccal(a)=\frac{\|a\|}{\mathrm{dist}(a,\Sigma)},
\end{equation*}
where $\| \; \|$ and $\mathrm{dist}$ are the norm and distance
induced by $\langle \ ,\ \rangle$. Note that for $\lambda\neq 0$
we have $\Ccal(\lambda a)=\Ccal(a)$. We can therefore work with
the $n$-dimensional real projective space $\PP^n$ as ambient
space. If we also denote by $\Sigma\subset \PP^n$ the image of the
ill-posed cone in projective space, then for $a\in \PP^n$ it
follows that
\begin{equation*}
  \Ccal(a)=\frac{1}{d_{\PP}(a,\Sigma)},
\end{equation*}
where $d_\PP(x,y)=\sin \alpha$, denotes the projective distance
between $x,y\in \PP^{n}$ ($\alpha$ being the angle between $x$ and
$y$).

The two-fold covering $p\colon S^n\rightarrow \PP^n$ induces a
measure $\nu$ on $\PP^n$ by means of
$\nu(B):=\frac{1}{2}\Vol_n(p^{-1}(B))$ for $B\subseteq \PP^n$,
where $\Vol_n$ is the $n$-dimensional volume on the sphere. Thus
$\nu(\PP^n)=\Ocal_n/2$, where $\Ocal_{n}:=\Vol_n
(S^n)=\frac{2\pi^{\frac{n+1}{2}}}{\Gamma (\frac{n+1}{2})}$.

For $0< \sigma \leq 1$ we denote by $B_\PP(a,\sigma)$ the open
ball of projective radius $\sigma$ around $a\in \PP^n$. It is
known that
\begin{equation*}
\nu(\B_\PP(a,\sigma))=\Ocal_{n-1}\cdot \I_{n}(\sigma),
\end{equation*}
where
\begin{equation}\label{jp}
   \I_{n}(\sigma):=\int_{0}^{\sigma}\frac{r^{n-1}}
     {\sqrt{1-r^2}}\;dr.
\end{equation}
The following bounds will prove useful on several occasions:
\begin{equation}\label{bound:ip}
  \frac{\sigma^n}{n}\leq \I_n(\sigma) \leq \min\left\{
  \frac{1}{\sqrt{1-\sigma^2}},\sqrt{\frac{\pi n}{2}}\right\}
  \cdot \frac{\sigma^n}{n}.
\end{equation}
For $a\in\PP^n$ and $\sigma\in (0,1]$ the uniform
measure on $\B_{\PP}(a,\sigma)$ is defined by 
\begin{equation}\label{uniform measure}
\nu_{a,\sigma}(B)=\frac{\nu(B\cap\B_{\PP}(a,\sigma))}
{\nu(\B_{\PP}(a,\sigma))}
\end{equation}
for all Borel-measurable $B\subseteq\PP^n$.

\subsection{Uniform smoothed analysis}

A reformulation of the main result in~\cite{BuCuLo:07} in the
projective space setting can be written as follows.\\

\begin{theorem}\label{mainthmold}
Let $\Ccal$ be a conic condition number with set of ill-posed
inputs $\Sigma\subset \PP^n$.  Assume that $\Sigma$ is contained
in the zero set in $\PP^n$ of homogeneous polynomials of degree at
most~$d$. Then, for all $\sigma\in(0,1]$ and all $t\geq
t_0=(2d+1)\frac{n}{\sigma}$,
\begin{equation*}
   \sup_{a\in \PP^n} \Prob_{z\in B_\PP(a,\sigma)}\{\Ccal(z)\geq
   t\}
     \leq  13\, d n\, \frac1{\sigma t} .
\end{equation*}
and
\begin{equation*}
   \sup_{a\in \PP^n}\,\Expect_{z\in B_\PP(a,\sigma)}[\ln\Ccal(z)]
      \leq 2\ln n+2\ln d +2\ln\frac1\sigma + 5,
\end{equation*}
where $\Prob$ and $\Expect$ are taken with respect to 
$\nu_{a,\sigma}$. 
\end{theorem}

As a consequence of this result, uniform smoothed analysis
results for the condition numbers of a variety of problems are obtained,
including linear equation solving, Moore-Penrose inversion,
eigenvalue computation and polynomial system solving. The bounds
obtained are consistently of the same order of magnitude as the
best bounds obtained previously by ad-hoc methods.

\subsection{Uniformly Absolutely Continuous Distributions}

In \cite{HM:06} a general boosting mechanism was developed that
allows extending any probabilistic analysis of a condition number with respect to
some chosen probability distribution over the input data to a more
general class of distributions. 

Let $\mu$ be a $\nu_{a,\sigma}$-absolutely continuous probability measure. 
Using the convention $\ln(0):=-\infty$ we define, for $\delta\in(0,1)$,
\begin{equation*}
\inf(\delta):=\inf\left\{\frac{\ln\mu(B)}{\ln\nu_{a,\sigma}(B)}: \,B\text{ is
}\text{Borel-measurable and }0<\nu_{a,\sigma}(B)\leq\delta\right\}
\end{equation*}
With these conventions, Theorem 2.2 of
\cite{HM:06} shows that
\begin{equation}\label{theorem 2.2}
\alpha_{\nu_{a,\sigma}}(\mu):=\lim_{\delta\rightarrow 0}\inf(\delta)\in[0,1].
\end{equation}
Absolute continuity alone ensures that all $\nu_{a,\sigma}$-null-sets must be
$\mu$-null-sets, but this does not imply that $\mu(B)$ is small
when $\nu_{a,\sigma}(B)$ is small and strictly positive. In contrast, when
$\alpha_{\nu_{a,\sigma}}(\mu)>0$ then \eqref{theorem 2.2} gives uniform upper
bounds on $\mu(B)$ in terms of $\nu_{a,\sigma}(B)$. Furthermore, the smaller
$\alpha$ gets, the larger the variation of $\mu$ in terms of
$\nu_{a,\sigma}$. If $\mu$ is $\nu_{a,\sigma}$-absolutely continuous and
$\alpha_{\nu_{a,\sigma}}(\mu)>0$, we therefore say that $\mu$ is {\em
uniformly} $\nu_{a,\sigma}$-absolutely continuous and call
$\alpha_{\nu_{a,\sigma}}(\mu)$ the {\em smoothness parameter} of $\mu$ with
respect to $\nu_{a,\sigma}$.

The following result, which easily follows from (\ref{theorem
2.2}), can be used to boost bounds on tail probabilities with
respect to $\nu_{a,\sigma}$ (as those in Theorem~\ref{mainthmold}) to obtain
similar bounds on
any uniformly $\nu_{a,\sigma}$-absolutely continuous probability measure $\mu$.\\

\begin{proposition}\label{eq:alphaup}
$\alpha_{\nu_{a,\sigma}}(\mu)$ is the largest nonnegative real number $\alpha$
for which it is true that
for all $\varepsilon > 0$ there exists $\delta_\varepsilon > 0$ such
that $\nu_{a,\sigma}(B)\leq\delta_{\varepsilon}$ implies  $\mu(B)\leq\nu_{a,\sigma}(B)^{\alpha
-\varepsilon}$.
\end{proposition}

\section{Smoothed analysis for adversarial distributions}
In this section we present our main result, namely an extension of
Theorem~\ref{mainthmold} to the case where we have a radially
symmetric distribution whose density has a pole at the point being
perturbed. We begin by introducing some notation.

Let $a\in\PP^n$ and $\sigma\in (0,1]$, and let 
$\nu_{a,\sigma}$ be the uniform measure on $\B_{\PP}(a,\sigma)$, 
as defined in \eqref{uniform measure}. 
Let $\mu$ be a $\nu_{a,\sigma}$-absolutely continuous
probability measure on $\PP^n$ with density $f(x)$. In other words, 
\begin{equation*}
  \mu(B)=\int_B f(x) \; \nu_{a,\sigma}(dx)
\end{equation*}
for all events $B$. 
Assume further that $f\colon \PP^n\rightarrow [0,\infty]$ is of
the form $f(x)=g(d_{\PP}(x,a))$, with a monotonically decreasing
function $g\colon [0,\sigma] \rightarrow [0,\infty]$ of the form
\begin{equation*}
g(r)=C_{\beta,\sigma}\cdot r^{-\beta}\cdot h(r),
\end{equation*}
with $\beta<n$, where
$C_{\beta,\sigma}=\I_n(\sigma)/\I_{n-\beta}(\sigma)$ and $h\colon
[0,\sigma] \rightarrow \R_{+}$ is a continuous function satisfying
$h(0)\neq 0$ and 
\begin{equation*}
\int_{0}^\sigma h(r)
\frac{r^{n-\beta-1}}{\sqrt{1-r^2}} \; dr=\I_{n-\beta}(\sigma),
\end{equation*}
so that $\mu$ is a probability measure on $B_{\PP}(a,\sigma)$. 
In other words, $f$ is radially symmetric around $a$
with respect to $d_{\PP}$ and has a pole of order $-\beta$ at $0$
in case $\beta>0$. The normalizing factor $C_{\beta,\sigma}$ is
chosen to make $h(r)=1$ a valid choice. 
Set $H:=\sup_{0\leq r\leq\sigma}h(r)$. Note that
$H\geq 1$, and that $H=1$ implies $h\equiv 1$.

It will be important to have expressions for 
$\nu_{a,\sigma}(B)$ and $\mu(B)$ when $B=B_{\PP}(a,\rho)$ is a
projective ball. In this situation we have
\begin{align}\label{eq:muradial}
  \mu(B_{\PP}(a,\rho))&=\frac{1}{\nu(B_{\PP}(a,\sigma))}
  \int_{B_{\PP}(a,\rho)}f(x)\; \nu(dx)\notag\\
  &=\frac{1}{\Ocal_{n-1}I_n(\sigma)}\cdot C_{\beta,\sigma}\cdot
  \Ocal_{n-1} \int_{0}^\rho r^{-\beta} h(r) \frac{r^{n-1}}{\sqrt{1-r^2}}\;
  dr\notag\\
 &=\frac{1}{\I_{n-\beta}(\sigma)}\int_{0}^\rho
  h(r)\frac{r^{n-\beta-1}}{\sqrt{1-r^2}} \; dr \\
  &\leq\left(\sup_{0\leq r\leq \rho}h(r)\right)\cdot\frac{\I_{n-\beta}(\rho)}
  {\I_{n-\beta}(\sigma)}.\notag
\end{align}
Similarly,
\begin{equation*}
 \mu(B_{\PP}(a,\rho))\geq
\left(\inf_{0\leq r\leq \rho} h(r)\right)\cdot
\frac{\I_{n-\beta}(\rho)}{\I_{n-\beta}(\sigma)}.
\end{equation*}
In particular,
\begin{equation}\label{eq:muradial2}
\nu_{a,\sigma}(B_{\PP}(a,\rho))=\frac{\I_{n}(\rho)}{\I_{n}(\sigma)}.
\end{equation}

The main result of this note is the following.\\

\begin{theorem}\label{thm:main}
Let $\Ccal$ be a conic condition number with set of ill-posed inputs
$\Sigma\subseteq \PP^n$, and assume $\Sigma$ is contained in a
projective hypersurface of degree at most $d$. Then
\begin{equation*}
  \Expect_{\mu}[\ln \Ccal]\leq
  2\ln(n)+\ln(d)+\ln\left(\frac{1}{\sigma}\right) +
  \ln\left(\frac{13\pi}{2}\right)+\frac{1}{1-\frac{\beta}{n}} \left(\ln\frac{2eH^2n}{\ln(\pi n/2)}\right).
\end{equation*}
\end{theorem}

This result applies to the variety of problems mentioned after
Theorem \ref{mainthmold}. The statement of the Theorem follows
from calculating the smoothness parameter $\alpha_\nu(\mu)$ and
the constants in Proposition~\ref{eq:alphaup}. These are given by
the following two lemmas, to be proven later.\\

\begin{lemma}\label{le:smoothness}
The smoothness parameter of $\mu$ with respect to $\nu_{a,\sigma}$
is given by $\alpha_{\nu_{a,\sigma}}(\mu)=1-\beta/n$.
\end{lemma}
\hspace{1cm}\\

For the statement of the next Lemma, let
$\varepsilon\in(0,1-\beta/n)$, and let
\begin{equation*}
  \rho_{\varepsilon}:= \sigma \cdot \left(\frac{1}{H}\cdot
  \sqrt{1-\left(\frac{2}{\pi n}\right)^{(1-\frac{\beta}{n}-\e)/(n\e)}}
    \right)^{\frac{1}{\e n}}\left(\sqrt{\frac{2}{\pi n}}\right)^{(1-\frac{\beta}{n}-\e)\frac{1}{\e n}}.
\end{equation*}
Set $\delta_{\varepsilon}:=\I_n(\rho_{\e})/\I_n(\sigma)$.\\

\begin{lemma}\label{le:boosting}
Let $B\subseteq \PP^n$ be such that $\nu_{a,\sigma}(B)\leq
\delta_{\e}$. Then $\mu(B)\leq
(\nu_{a,\sigma}(B))^{1-\frac{\beta}{n}-\e}$.
\end{lemma}
\hspace{1cm}\\

We are now ready to prove the main result.

\Proofof{Theorem \ref{thm:main}} Setting
$\e=\frac{1}{2}(1-\frac{\beta}{n})$ and using the bounds
(\ref{bound:ip}) we obtain
\begin{equation}\label{bound:deltae}
  \frac{2}{\pi n}\left(\frac{1}{H}\cdot \sqrt{1-\left(\frac{2}{\pi n}\right)^{\frac{1}{n}}}\right)^{\frac{2}{1-\frac{\beta}{n}}}\leq \delta_\e\leq
  \left(\frac{1}{H}\cdot \sqrt{1-\left(\frac{2}{\pi n}\right)^{\frac{1}{n}}}\right)^{\frac{2}{1-\frac{\beta}{n}}}.
\end{equation}

From Theorem~\ref{mainthmold} it follows that for all $t\geq
t_0:=\ln[(1+2d)n/\sigma]$,
\begin{equation}\label{BCL}
\Prob_{\nu_{a,\sigma}}\{\ln\Ccal>t\}\leq
\frac{13dn}{\sigma}e^{-t}.
\end{equation}
Set
\begin{equation*}
t_{\varepsilon}:=\ln\left(\frac{13dn}{\sigma \cdot
\delta_\e}\right)= \ln\left(\frac{13dn}{\sigma}\right) +\ln
(\delta_\e^{-1}).
\end{equation*}
Using (\ref{bound:deltae}) we obtain
\begin{equation*}
  \ln\left(13\frac{dn}{\sigma}\right)\leq
  t_{\e}-\frac{2}{1-\frac{\beta}{n}}\ln
  \left(\frac{H}{\sqrt{1-\left(\frac{2}{\pi n}\right)^{
\frac{1}{n}}}}\right)\leq \ln\left(13\frac{\pi}{2}\frac{d
    n^2}{\sigma}\right).
\end{equation*}
The lower bound shows that $t_{\varepsilon}>t_0$, so that for all
$t\geq t_{\varepsilon}$,
\begin{equation*}
\nu_{a,\sigma}\left(\left\{x:\,\ln\Ccal(x)>t\right\}\right)
=\Prob_{\nu_{a,\sigma}}\{\ln\Ccal>t\}\leq
\frac{13dn}{\sigma}e^{-t} \leq\delta_{\varepsilon}.
\end{equation*}
Applying Lemma~\ref{le:boosting}, it follows that for $t\geq
t_{\varepsilon}$,
\begin{equation*}
\Prob_{\mu}\{\ln\Ccal>t\}=\mu\left(\left\{x:\,\ln\Ccal(x)>t\right\}\right)
\leq\left(\frac{13dn}{\sigma}e^{-t}
\right)^{\frac{1}{2}(1-\frac{\beta}{n})},
\end{equation*}
and hence,
\begin{align*}
\Expect_{\mu}[\ln\Ccal]&=\int_{0}^{\infty}\Prob_{\mu}\{\ln\Ccal>t\}dt\\
&\leq\int_{0}^{t_{\varepsilon}}1\,dt+\int_{t_{\varepsilon}}^{\infty}
\left(\frac{13dn}{\sigma}e^{-t}
\right)^{\frac{1}{2}(1-\frac{\beta}{n})}\,dt\\
&= t_\varepsilon + \frac{2
\delta_{\varepsilon}^{\frac{1}{2}(1-\frac{\beta}{n})}}{1-\frac{\beta}{n}}.
\end{align*}
Using the bounds on $t_\e$ and $\delta_\e$ we get
\begin{equation*}
\Expect_{\mu}[\ln \Ccal]\leq
2\ln(n)+\ln(d)+\ln\left(\frac{1}{\sigma}\right)+\ln\left(\frac{13\pi}{2}\right)+\frac{2}{1-\frac{\beta}{n}}\left(\ln\left(\frac{H}{\sqrt{1-\left(\frac{2}{\pi
      n}\right)^{\frac{1}{n}}}}\right)+ \frac{\sqrt{1-\left(\frac{2}{\pi
      n}\right)^{\frac{1}{n}}}}{H}\right).
\end{equation*}
A small calculation shows that $\left(1-\left(\frac{2}{\pi
      n}\right)^{\frac{1}{n}}\right)^{-1/2}\leq \sqrt{\frac{2n}{\ln(\pi n/2)}}$. This completes the proof.
\endProofof

\subsection{Proofs of Lemmas \ref{le:smoothness} and \ref{le:boosting}}\label{se:tailbounds}
The content of the following Lemma,
needed for calculating the smoothness parameter, should be intuitively clear.\\

\begin{lemma}\label{le:among}
Let $0<\delta<1$. Then among all measurable sets
$B\subseteq\B_{\PP}(a,\sigma)$ with
$0<\nu_{a,\sigma}(B)\leq\delta$, $\mu(B)$ is maximized by
$\B_{\PP}(a,\rho)$ where $\rho\in(0,\sigma)$ is chosen so that
$\nu_{a,\sigma}(\B_{\PP}(a,\rho))=\delta$.
\end{lemma}

\begin{proof}
It clearly suffices to show that
\begin{equation*}
\int_{B}f(x)\;\nu_{a,\sigma}(dx)
\leq\int_{\B_{\PP}(a,\rho)}f(x)\;
\nu_{a,\sigma}(dx)
\end{equation*}
for all Borel sets $B\subset\B_{\PP}(a,\sigma)$ such that
$\nu_{a,\sigma}(B)=\delta$. Indeed, we have
\begin{align}
\int_{B}f(x)&\;\nu_{a,\sigma}(dx)=
\int_{B\cap \B_{\PP}(a,\rho)}f(x)\;\nu_{a,\sigma}(dx)
+\int_{B\setminus\B_{\PP}(a,\rho)}f(x)\;\nu_{a,\sigma}(dx)
\nonumber\\
&\leq\int_{B\cap \B_{\PP}(a,\rho)}f(x)\;\nu_{a,\sigma}(dx)
+g(\rho)\;\nu_{a,\sigma}(B\setminus\B_{\PP}(a,\rho))\nonumber\\
&=\int_{B\cap \B_{\PP}(a,\rho)}f(x)\;\nu_{a,\sigma}(dx)
+g(\rho)\;\nu_{a,\sigma}(\B_{\PP}(a,\rho)\setminus B)\label{delta}\\
&\leq\int_{B\cap \B_{\PP}(a,\rho)}f(x)\;\nu_{a,\sigma}(dx)
+\int_{\B_{\PP}(a,\rho)\setminus B}f(x)\;\nu_{a,\sigma}(dx)
\nonumber\\
&=\int_{\B_{\PP}(a,\rho)}f(x)\;\nu_{a,\sigma}(dx),\nonumber
\end{align}
where we have used $\nu_{a,\sigma}(\B_{\PP}(a,\rho))=\delta=\nu_{a,\sigma}(B)$
in \eqref{delta}. This proves our claim.
\end{proof}

Even though $\rho$ is a function of $\delta$, we will not reflect
this notationally in the sequel.

\Proofof{Lemma \ref{le:smoothness}} From (\ref{eq:muradial}),
(\ref{eq:muradial2}) and (\ref{bound:ip}) we get the bounds of the form 
\begin{align}
\frac{1}{C_1}\cdot \rho^n
&\leq\nu_{a,\sigma}(\B_{\PP}(a,\rho))
\leq C_1\cdot \rho^n,\label{bounds on delta}\\
\inf_{0\leq r\leq \rho} h(r) \cdot \frac{1}{C_2} \cdot \rho^{n-\beta}
&\leq\mu(\B_{\PP}(a,\rho))\leq \sup_{0\leq r\leq \rho}h(r) \cdot C_2\cdot \rho^{n-\beta},\label{bounds on mu}
\end{align}
where the constants $C_i$ do not depend on $\rho$.

We thus have (using Lemma ~\ref{le:among})
\begin{align*}
\alpha_{\nu_{a,\sigma}}(\mu)
&=\lim_{\delta\rightarrow 0}\inf\left\{\frac{\ln\mu(B)}{\ln\nu_{a,
\sigma}(B)}:\,B\text{ measurable},\,0<\nu_{a,\sigma}(B)\leq\delta\right\}\\
&=\lim_{\rho\rightarrow 0}\frac{\ln\mu(\B_{\PP}(a,\rho))}{\ln\nu_{a,\sigma}
(\B_{\PP}(a,\rho))}\\
&\begin{cases}\leq\lim_{\rho\rightarrow 0}\frac{\ln (\inf h(r)/C_2)+(n-\beta)\ln \rho}
{\ln \left(C_1\right)+n\ln \rho}
=1-\frac{\beta}{n}\\[6pt]
\geq\lim_{\rho\rightarrow 0}\frac{\ln (C_2
\cdot \sup h(r))+(n-\beta)\ln \rho}{-\ln C_1+n\ln \rho}
=1-\frac{\beta}{n}.
\end{cases}
\end{align*}
This concludes the proof.
\endProofof

\Proofof{Lemma \ref{le:boosting}} Since sets of the form
$\B_{\PP}(a,\rho)$ maximise $\mu(B)$ among all measurable sets
$B\subseteq\B_{\PP}(a,\sigma)$ such that $\nu_{a,\sigma}(B)
\leq\delta$ for any $\delta$, we may w.l.o.g. assume
$B=B_{\PP}(a,\rho)$. By (\ref{eq:muradial}) and
(\ref{eq:muradial2}) our task amounts to showing
\begin{equation*}
H \cdot \frac{\I_{n-\beta}(\rho)}{\I_{n-\beta}(\sigma)}\leq \left(\frac{\I_n(\rho)}{\I_n(\sigma)}\right)^{1-\frac{\beta}{n}-\e}
\end{equation*}
for $\rho\leq \rho_{\e}$. And indeed, using the bounds
(\ref{bound:ip}), we get
\begin{align*}
  H \cdot \frac{\I_{n-\beta}(\rho)}{\I_{n-\beta}(\sigma)} &\leq
  H \frac{1}{\sqrt{1-\rho^2}}\cdot \left(\frac{\rho}{\sigma}\right)^{n-\beta}\\
  &\leq
  H\frac{1}{\sqrt{1-\rho^2}}\cdot \left(\left(\frac{\rho}{\sigma}\right)^{n}\right)^{1-\frac{\beta}{n}-\e}\left(\frac{\rho_\e}{\sigma}\right)^{\e
    n}\\
&\leq \frac{\sqrt{1-\left(\frac{2}{\pi n}\right)^{(1-\frac{\beta}{n}-\e)/(n\e)}}}{\sqrt{1-\rho^2}}\cdot
  \left(\sqrt{\frac{2}{\pi n}}\left(\frac{\rho}{\sigma}\right)^{n}\right)^{1-\frac{\beta}{n}-\e}\\
&\leq \frac{\sqrt{1-\left(\frac{2}{\pi n}\right)^{(1-\frac{\beta}{n}-\e)/(n\e)}}}{\sqrt{1-\rho^2}}\cdot
  \left(\frac{\I_n(\rho)}{\I_n(\sigma)}\right)^{1-\frac{\beta}{n}-\e},
\end{align*}
where for the last inequality we use the bounds (\ref{bound:ip}) again.
Moreover, we have
\begin{equation*}
  \rho \leq \rho_{\e}\leq \left(\sqrt{\frac{2}{\pi
      n}}\right)^{(1-\frac{\beta}{n}-\e)\frac{1}{\e n}}.
\end{equation*}
Therefore, $\sqrt{1-\left(\frac{2}{\pi n}\right)^{(1-\frac{\beta}{n}-\e)\frac{1}{\e n}}}\leq \sqrt{1-\rho^2}$, completing
the proof.
\endProofof

\goodbreak
{\small

}
\end{document}